\documentclass[12pt]{article}
\usepackage{amsmath}
\usepackage{amssymb}
\newsymbol \blackbox 1004
\newcommand{\eh}{\hfill}\newlength{\sperr}

\newenvironment{proof}{{\settowidth{\sperr}{\bf\rm
Proof}%
\par\addvspace{0.3cm}\noindent\parbox[t]{1.3\sperr}
{\bf\rm P\eh r\eh o\eh o\eh f\eh }%
}}{\nopagebreak\mbox{}
$\blackbox$\par\addvspace{0.3cm}}

\def\a{\alpha}
\def\b{\beta}
\def\g{\gamma}
\def\G{\Gamma}

\def\s{\sigma}
\def\la{\lambda}
\def\om{\omega}
\def\Om{\Omega}
\def\S{\Sigma}

\def\t{\theta}

\def\vp{\varphi}
\def\ve{\varepsilon}
\def\wh{\widehat}
\def\wt{\widetilde}

\def\BC{{\mathbb C}}
\def\BR{{\mathbb R}}

\newtheorem{Pa}{Paper}[section]
\newtheorem{Tm}[Pa]{{\bf Theorem}}
\newtheorem{La}[Pa]{{\bf Lemma}}
\newtheorem{Cy}[Pa]{{\bf Corollary}}
\newtheorem{Rk}[Pa]{{\bf Remark}}

\newtheorem{Dn}[Pa]{{\bf Definition}}
\newtheorem{Pn}[Pa]{{\bf Proposition}}
\newcommand{\CC}
{{\mathchoice {\setbox0=\hbox{$\displaystyle\rm
C$}\hbox{\hbox
to0pt{\kern0.4\wd0\vrule height0.9\ht0\hss}\box0}}
{\setbox0=\hbox{$\textstyle\rm C$}\hbox{\hbox
to0pt{\kern0.4\wd0\vrule height0.9\ht0\hss}\box0}}
{\setbox0=\hbox{$\scriptstyle\rm C$}\hbox{\hbox
to0pt{\kern0.4\wd0\vrule height0.9\ht0\hss}\box0}}
{\setbox0=\hbox{$\scriptscriptstyle\rm C$}\hbox{\hbox
to0pt{\kern0.4\wd0\vrule height0.9\ht0\hss}\box0}}}}

\title{Semiseparable integral operators  and explicit solution
of an inverse problem for the  skew-self-adjoint
Dirac-type system}

\author{B. Fritzsche, B. Kirstein,
A.L. Sakhnovich}

\date{}
\begin{document}
\maketitle

\begin{abstract} Inverse problem to recover the skew-self-adjoint Dirac-type system from the
generalized Weyl matrix function is treated in the paper. Sufficient conditions under which  the unique solution
of the inverse problem exists, are formulated in terms of the Weyl function and a procedure to
solve the inverse problem is given. The case of the generalized Weyl functions of the form $\phi(\lambda)\exp\{-2i\lambda D\}$,
where $\phi$ is a strictly proper rational matrix function and $D=D^* \geq 0$ is a diagonal matrix,
is treated in greater detail.  Explicit formulas for the inversion of the corresponding semiseparable
integral operators and recovery of the Dirac-type system are obtained for this case.

\end{abstract}

{MSC(2000) 34A55; 34B20; 47G10; 34A05}


\section{Introduction} \label{intro}
\setcounter{equation}{0}
The  skew-self-adjoint Dirac-type system
\begin{equation}      \label{0.1}
\frac{d}{dx}u(x, \lambda )=\Big(i \lambda j+j
V(x)\Big)u(x,
\lambda ), \quad x \geq 0,
\end{equation}
where
\begin{equation}  \label{0.2'}
j = \left[
\begin{array}{cc}
I_{p} & 0 \\ 0 & -I_{p}
\end{array}
\right], \hspace{1em} V= \left[\begin{array}{cc}
0&v\\v^{*}&0\end{array}\right],
\end{equation}
$I_p$ is the $p \times p$ identity matrix, and $v$ is a $p \times
p$ matrix function, is actively studied in analysis and soliton theory
(see, for instance, \cite{APT, FT} and the references therein). System (\ref{0.1}) is an auxiliary linear system
for the focusing matrix NLSE, sine-Gordon and other important integrable equations.

The inverse problem to recover  a  self-adjoint Dirac type system from its
Weyl or spectral function is closely related to the inversion of the integral
operators  with difference kernels, see \cite{DI, Kre1, SaA3, SaL2, SaL3} and various
references. For the discrete analogues of Dirac systems, Toeplitz matrices appear
instead of the operators with difference kernels \cite{DGK2, Dy, FKRS, Si2}.
(Various results on Toeplitz matrices and related $j$-theory one can find, for instance, in \cite{BS,
DFK, FrKi, FrKi2}.)

When the Weyl functions of the self-adjoint Dirac type system are rational,
one can solve the inverse problem explicitly.  One of the approaches to
solve the inverse problem explicitly is connected with a version of the
B\"acklund-Darboux transformation and some notions from system theory
\cite{GKS98, GKS02}. (See also  \cite{FKRS, FKS1, KaS, MST} for this approach, and see \cite{GT}
and the references therein for explicit formulas for the radial Dirac equation.)
Another method is to apply the general theory.
It proves \cite{AGKLS} that for the case of rational Weyl functions the
corresponding operators with difference kernels can be inverted
explicitly by formulas from \cite{BGK}.

The case of the skew-self-adjoint Dirac  type system
with the rational Weyl function was treated in  \cite{GKS2} .
It was shown  that any strictly proper
rational $p \times p$ matrix  function is the Weyl function 
of  a  skew-self-adjoint Dirac  type system on semi-axis
and the solution of the inverse problem was constructed explicitly
similar to the self-adjoint case treated in \cite{GKS98}.

The analogues of  the operators with difference kernel for the 
skew-self-adjoint system (\ref{0.1}) 
are    bounded operators
$S_l$ in $L^2_p(0,l)$  ($0<l<\infty$), which have the 
form \cite{SaA2, SaA8}
\begin{equation}\label{0.2}
S_l f=Sf=f(x)+ \frac{1}{2}\int_0^l
\int^{x+t}_{|x-t|} k\Big(
\frac{r+x-t}{2} \Big) k \Big( \frac{r+t-x}{2} \Big)^*d r f(t) d
t,
\end{equation}
where 
$\sup_{0<x<l}\|k(x)\|<\infty$. The kernel of the operator $S_l$
is denoted by $K$:
\begin{equation}\label{v27}
K(x,t)= \frac{1}{2}
\int^{x+t}_{|x-t|} k\Big(
\frac{r+x-t}{2} \Big) k \Big( \frac{r+t-x}{2} \Big)^*d r .
\end{equation}

In this paper we show that for a Weyl function
of the form
\begin{equation}\label{a1}
\vp(\la)=\phi(\la)\exp\{-2i  \lambda D\}R, \quad D \geq 0,
\end{equation}
where $\phi$ is a strictly proper rational $p \times p$ matrix  function, 
$D$ is a $p \times p$ diagonal  matrix, and $R$ is a $p \times p$ 
unitary matrix,
the corresponding operator $S$ is semiseparable. 
Using  results on the inversion of the semiseparable operators,
the inverse problem
to recover the system from $\vp$ is solved explicitly.  Putting $D=0$,
we get the subcase of rational Weyl functions. 
Some definitions and results
for the general type (non-explicit) case of inverse problem including Theorem \ref{TmssaD}
and the important
formula   (\ref{1.8d3})  are also  new.
The semiseparable matrices and operators are actively studied
(see, for instance, \cite{EG02, GK84, GGK1, VVGM}), and their application to inverse problems
is of interest, too.

Various definitions and results on a general type inverse problem for 
the skew-self-adjoint Dirac  type system and on explicit solutions of
the inverse problem, when the Weyl functions are proper rational,
are given in Section \ref{Prel}. Some properties of the operator
$S_l$ of the form (\ref{0.2}) are studied in Section \ref{fact}.
The explicit solution of the inverse problem for the Weyl functions of the form
(\ref{a1})  is contained in Section \ref{ISP}.

We denote  by $\BR$ the real axis, by $\BR_+$
the positive semi-axis, by
$\BC$ the complex plane,
and by $\BC_+$  the open upper halfplane $\Im \la>0$.
The class of bounded linear operators
acting from $H_1$ into $H_2$ is denoted by    $\{H_1, \,  H_2\}$,
 the identity operators are denoted by $I$,
and spectrum is denoted by $\s$.

\section{Inverse problem. Preliminaries} \label{Prel}
\setcounter{equation}{0}
First, normalize the fundamental solution $u(x, \la)$ of system 
(\ref{0.1})  by the initial condition
\begin{equation} \label{v0}
u(0, \la)=I_{2p}.
\end{equation}
If
\begin{equation} \label{1.2}
\sup_{0<x<\infty}\|v(x)\| \leq M,
\end{equation}
the  unique
$p \times p$ Weyl matrix function $\vp(\lambda)$  of the  skew-self-adjoint Dirac type system
 (\ref{0.1})  on the semi-axis $[0, \, \infty)$ can be defined
\cite{SaA1} (see also 
\cite{CG2, GKS2, SaA8}) by the inequality
\begin{equation} \label{1.1}
\int_0^\infty \left[ \begin{array}{lr}    \varphi (\lambda)^* &
I_p
\end{array} \right]
  u(x, \lambda)^*
u(x, \lambda)
\left[ \begin{array}{c}
  \varphi (\lambda) \\ I_p \end{array} \right] dx < \infty,
\end{equation}
which holds for all $\lambda$ in the halfplane  $\Im \lambda < -M<0$.
Under condition (\ref{1.2}) such a Weyl function always exists.

Consider the case of the so called pseudo-exponential potentials \cite{GKS2},
which are denoted by the acronym PE.
A potential $v \in $PE is determined by  three parameter matrices, that is, by the
$n \times n$ matrix $\a$ ($n>0$) and  two $p \times n$ matrices
$\t_1$ and $\t_2$, which satisfy the identity
\begin{equation}\label{1.9}
\a -\a^*=i(\t_1\t_1^*+\t_2\t_2^*).
\end{equation}
The pseudo-exponential potential has the form
\begin{equation}\label{1.10}
v(x)=2 \theta _{1}^{*}e^{ix \alpha ^{*}}\S(x)^{-1}e^{ix \alpha }
\theta _{2}, 
\end{equation}
where
\begin{equation}\label{1.11}
\S(x)=I_n+\int_0^x \Lambda (t)j \Lambda (t)^*dt,  \quad  \Lambda (x)= \left[ \begin{array}{lr}
e^{-ix \alpha } \theta_{1} & e^{ix \alpha } \theta_{2}
\end{array}
\right].
\end{equation}
By Proposition 1.4 in \cite{GKS2} the pseudoexponential potential $v$, i.e., the potential
 given by  (\ref{1.10})  is bounded on the semi-axis.
The Weyl function of the system  (\ref{0.1}) with $v \in$PE
is a rational matrix function, which is also expressed in terms of the parameter
matrices \cite{GKS2}:
\begin{equation} \label{1.12}
\vp ( \lambda )=i \theta_{1}^{*}( \lambda I_{n} - \beta )^{-1} \theta_{2}, \quad \b=\a - i\t_2 \t_2^*.
\end{equation}
In spite of the requirement
\begin{equation} \label{v23}
 \b-\b^*=i(\t_1 \t_1^* - \t_2 \t_2^*),
\end{equation}
which is implied by   the equalities (\ref{1.9}) and  $\b=\a - i\t_2 \t_2^*$, any strictly
proper rational matrix function
can be presented in the form (\ref{1.12}).
The inverse problem to recover $v$ from the strictly 
proper rational matrix function $\vp$ is solved explicitly in \cite{GKS2},
using a minimal realization  of $\vp$ and formula (\ref{1.10}).

When  (\ref{1.2}) is true, inequality (\ref{1.1}) implies other inequalities:
\begin{equation}  \label{v1}
\sup_{x \leq l, \, \Im \la<-M} \left\|e^{i x \la} u(x, \la)
\left[\begin{array}{c}
\vp(\la) \\ I_{p}
\end{array}
\right]\right\|<\infty \quad {\mathrm{for}} \, {\mathrm{all}} \, \, 0<l< \infty ,
\end{equation}
which can be treated as a more general definition of the Weyl function.
\begin{Dn} \label{Dn1} 
Let the system (\ref{0.1}) be given on the semi-axis $[0, \, \infty)$.
Then a $p \times p$ matrix function $\vp(\lambda)$ analytic  in
some halfplane $\Im \lambda < -M<0$ is called a Weyl function of
this system, if inequalities (\ref{v1}) hold.
\end{Dn}
If
\begin{equation} \label{v2}
\sup_{0<x<l}\|v(x)\| <\infty  \quad {\mathrm{for}} \, {\mathrm{all}} \, \, 0<l< \infty,
\end{equation}
then there is at most one Weyl function.
\begin{Dn} \label{Dn2} 
The inverse spectral problem (ISP)  for system 
(\ref{0.1})  on the semi-axis is the problem to recover $v(x)$
satisfying (\ref{v1}) and  (\ref{v2}) from the Weyl function $\vp$.
\end{Dn}
For an analytic matrix function $\vp$ satisfying the condition
\begin{equation} \label{v3}
\sup_{\Im \la <-M}\|\la^2\big(\vp(\la)-\a/\la \big)\| <\infty , 
\end{equation}
where $\alpha$ is some $p \times p$ matrix, the solution of the inverse problem
always exists (see Lemma 1  \cite{SaA2} for the scalar version of this result
and the matrix case can be proved quite similar).

The general (non-explicit) procedure to solve ISP is described in \cite{SaA0, SaA1, SaA2, SaA8}.
Fix a positive value $l$ ($0<l< \infty$).
The first step  to solve ISP is to recover
a $p \times p$ matrix function $s(x)$  with
the entries from $L^2(0,l)$ ($l<\infty$), i.e., $s(x) \in L^2_{p \times p}(0,l)$ via the
Fourier transform. That is, we put
\begin{equation}\label{1.3}
\displaystyle{s(x)= \frac{i}{2 \pi}e^{- \eta x}{\mathrm{
l.i.m.}}_{a \to \infty} \int_{- a}^{a}e^{i \xi
x} \lambda^{-1} \vp(\lambda /2) d \xi \quad (\lambda= \xi +i \eta , \quad \eta
<-2M),}
\end{equation}
the limit l.i.m. being the limit in $L^2(0,l)$. As (\ref{1.3}) has sense for any $l<\infty$
the matrix function $s(x)$ is defined on the non-negative real  semi-axis $x\geq 0$.  Moreover,
it is easily checked that $s$ is absolutely continuous, it does not depend on the choice of $\eta <-
2M$, $s^{\prime}$ is bounded on any finite interval, and $s(0)=0$.
To define the operator $S_l$ we substitute $k(x)=s^{\prime}(x)$
into  (\ref{0.2}).

Next, denote the $p\times 2p$ block rows of $u$ by $\om_1$ and $\om_2$:
\begin{equation}\label{1.5}
\om_1(x)=[ I_p \quad 0]u(x, 0), \quad  \om_2(x)=[0 \quad I_p]u(x, 0).
\end{equation}
It follows from (\ref{0.1}) that $u(x,0)^*u(x,0)=I_{2p}$. Hence, by (\ref{0.1})  and  (\ref{1.5}) we have
\begin{equation}\label{1.6}
v(x)=\om_1^{\prime}(x) \om_2(x)^*,
\end{equation}
and $\om_1$, $\om_2$ satisfy the equalities
\begin{equation} \label{1.7}
\om_1(0)=[I_p \quad 0], \quad \om_1 \om_1^* \equiv I_p, \quad
\om_1^{\prime}
\om_1^* \equiv 0, \quad 
\om_1 \om_2^* \equiv 0.
\end{equation}
It is immediate that $\om_1$ is uniquely recovered from $\om_2$ using  (\ref{1.7}).

Finally, we obtain $\om_2$ via the formula
\begin{equation}\label{1.8}
\om_2(l)=[0 \quad I_p]- \int_0^l \Big(S_l^{-1} s'(x)\Big)^*[I_p
\quad s(x)]dx \quad (0<l< \infty),
\end{equation}
where $S_l^{-1}$ is applied to $ s^{\prime}$ columnwise.

From the considerations in \cite{SaA1, SaA2} (see also similar constructions in \cite{SaAnw},
where the Weyl theory for the linear system auxiliary to the nonlinear optics equation
is treated) it follows that one can solve ISP under requirements on $\vp$
and $s(x)$ weaker than (\ref{v3}). Namely, we assume
\begin{equation}\label{v4}
\sup_{\Im \la < -M}\|\vp(\la)\| <\infty,
\end{equation}
\begin{equation}\label{v5}
\vp(\la) \in L^2_{p \times p}(- \infty, \, \infty), \quad  \la= \xi + i \eta \,\,  (- \infty <\xi< \infty )
\, \, {\mathrm{for}} \, {\mathrm{all}} \, \, \eta < - M,
\end{equation}
\begin{equation}\label{v6}
s(0)=0, \quad \sup_{0<x<l}\|k(x)\|<\infty 
\quad {\mathrm{for}} \, {\mathrm{all}} \, \, 0<l< \infty, \quad
k(x):=s^{\prime}(x),
\end{equation}
\begin{equation}\label{v7}
\int_0^{\infty}e^{-cx}\| k(x) \| dx < \infty
\end{equation}
for some $c>0$.
\begin{Tm}\label{TmssaD} 
Let  the matrix function $\vp$ be analytic in the halfplane $\Im \la < -M$ and   satisfy
the relations  (\ref{v4}) and (\ref{v5}).
 Let also the  matrix function $s(x)$   defined via $\vp$ by formula  (\ref{1.3}) 
 be absolutely continuous and satisfy (\ref{v6}) and  (\ref{v7}).
Then ISP  has a unique solution,
which is given by formulas  (\ref{1.6})- (\ref{1.8}), where $S_l \geq I$ has the form  (\ref{0.2})
with $k=s^{\prime}$.
\end{Tm}
%%%%%%%%%%%%%%%%%%%%%%%%%%%%%
\section{Factorization of $S$ and operator   identity} \label{fact}
\setcounter{equation}{0}
Consider again the operator $S=S_l$. It is easy to see that functions, which are bounded on
the interval, can be approximated in the $L^1$-norm by the continuous functions.
As $k=s'$ is bounded on the finite intervals, one can see that the kernel $K$ of $S$, which is given 
by (\ref{v27}), is continuous with respect to $x$ and $t$. Hence, the kernel
of $S_l^{-1}$ is continuous with respect to $x$, $t$, and $l$ (\cite{GoKr}, p. 185).
Therefore, $S_l^{-1}k$ has the form $\Big(S_l^{-1}k\Big)(x)=k(x)+k_1(x)$, where
$k_1$ is continuous, and  the matrix function $\Big(S_l^{-1}k\Big)(l)$ is well-defined:
\begin{equation}\label{3v1}
\Big(S_l^{-1}k\Big)(l)=k(l)+k_1(l)=k(l)+ \lim_{x \to l-0}\Big(\big(S_l^{-1}k\big)(x)-k(x)\Big) \quad (0<l<\infty).
\end{equation}

To express $v$ in terms of  $\Big(S_l^{-1}k\Big)(l)$ we need some preparations.
According to \cite{SaL0}  there are triangular operators  $\wh V_l \in \{L^2_p(0,l), \, L^2_p(0,l)\}$,  such that
\begin{equation}\label{3v2}
(\wh V_lf)(x)=f(x)+\int_0^x\wh V_-(x,t) f(t)dt,  \quad
\wh  V_l A \wh  V_l^{-1}=i \om_1(x)\int_0^x\om_1(t)^* \, \cdot \, dt,
\end{equation}
$\wh V_-(x,t)$ does not depend on $l$, and the operators $\wh V_l$ and $\wh V_l^{-1}$ map functions with 
bounded derivatives into functions with bounded derivatives. Moreover,
as bounded functions on an interval can be approximated in the  $L^1$-norm by the continuous functions,
it follows from the construction in \cite{SaL0}  that $\wh V(x,t)$ ($x \geq t$) is continuous with respect to $x$ and $t$.

Next, introduce the operator
\begin{equation}\label{3v3}
(\wt V_lf)(x)=f(x)+\int_0^x\wt V_-(x-t) f(t)dt,  \quad  \wt V_-(x):=\frac{d}{d x}\Big(\wh V_l^{-1}\om_{11}\Big)(x),
\end{equation}
where $\om_{11}$ is the first $p \times p$ block of $\om_1$,
and put 
\begin{equation}\label{3v4}
V_l:=\wh  V_l \wt V_l=I+\int_0^x  V_-(x,t) \, \cdot \, dt.
\end{equation}
 It is easy to see that $\wt V_l A=A\wt V_l$, and so the second equality in (\ref{3v2})
yields
\begin{equation}\label{3v5}
  V_l A   V_l^{-1}=i \om_1(x)\int_0^x\om_1(t)^* \, \cdot \, dt.
\end{equation}
By (\ref{1.5}) we see that $\om_{11}(0)=I_p$. Hence, using definition  (\ref{3v3}) one gets
\begin{eqnarray}\label{3v6}
\wt V_l I_p&=&I_p+\int_0^x \Big(\frac{d}{d x}\big(\wh V_l^{-1}\om_{11}\big)\Big)(x-t)dt 
\\ \nonumber
&=&I_p+\int_0^x \Big(\frac{d}{d x}\big(\wh V_l^{-1}\om_{11}\big)\Big)(t)dt=\big(\wh V_l^{-1}\om_{11}\big)(x).
\end{eqnarray}
Formula (\ref{3v6}) implies $V_l^{-1}\om_{11}=I_p$.
Moreover, from   \cite{SaA1, SaA8} it follows that under the conditions of Theorem \ref{TmssaD} the equalities
\begin{equation}\label{1.8d0'}
 \big(V_l^{-1}\om_1\big)(x)=[I_p \quad s(x)]
\end{equation}
and
\begin{equation}\label{1.8d3'}
S_l^{-1}=V_l^*V_l
\end{equation}
are also true.
\begin{Rk}\label{RkV}  Under the conditions of Theorem \ref{TmssaD}  we have
\begin{equation}\label{1.8d3}
v(l)=\big(S_l^{-1}s^{\prime}\big)(l).
\end{equation}
Indeed, using  (\ref{1.8d3'}) and
changing variables $l$ and $x$  into
$x$ and $t$, correspondingly, we rewrite  (\ref{1.8}) in the form
\begin{equation}\label{1.8d4}
\om_2(x)=[0 \quad I_p]- \int_0^x\Big(V_x s'\Big)(t)^*V_x[I_p
\quad s(t)]dt.
\end{equation}
As $V_-$ does not depend on $l$ we have $\big(V_x s'\big)(t)=\big(V_l s'\big)(t)$ for $t\leq x \leq l $.
Thus, according to   (\ref{1.8d0'}) and   (\ref{1.8d4}), we get
\begin{equation}\label{1.8d5}
\om_2^{\prime}(x)=-  \big(V_l s'\big)(x)^*\om_1(x).
\end{equation}
Multiplying both sides of (\ref{1.8d5}) by $\om_1^*$
from the right and taking into account (\ref{1.6}) and (\ref{1.7}),  one derives
$-v(x)^*=-  \big(V_l s'\big)(x)^*$, i.e., the equality
\begin{equation}\label{1.8d6}
 v(x)=\big(V_ls^{\prime}\big)(x)
\end{equation}
is true.  As $\wh V(x,t)$ is continuous, taking into account (\ref{3v3}) and (\ref{3v4})
we see that $\big(V_ls^{\prime}\big)(x)-s^{\prime}(x)$ is continuous.
 It is also immediate from  (\ref{3v4}) that 
\begin{equation}\label{1.8d1'}
(V_l^*f)(x)=f(x)+\int_x^lV_-(t,x)^* f(t)dt.
\end{equation}
Hence, according to (\ref{3v1}), (\ref{1.8d3'}), and  (\ref{1.8d1'}) we have
\begin{equation}\label{3v7}
\Big(S_l^{-1}s^{\prime}\Big)(l)=\Big(V_l 
s^{\prime}\Big)(l).
\end{equation}
Finally, formula  (\ref{1.8d3}) follows from (\ref{1.8d6}) and  (\ref{3v7}).
\end{Rk}

By  (\ref{1.8d3'})  the equality
\[
AS-SA^*=V_l^{-1} \Big( V_lAV_l^{-1} -\big(  V_lAV_l^{-1}\big)^*\Big)\big(V_l^{-1} \big)^*
\] 
is valid for $S=S_l$. Therefore,  taking into account (\ref{3v5}) and (\ref{1.8d0'}) one can see that
 $S$ satisfies the operator identity
\begin{equation}\label{0.4}
AS-SA^*=i \Pi \Pi^*, \quad  \Pi=[\Phi_1 \quad \Phi_2], \quad 
\Phi_1 g \equiv g, \quad \Phi_2 g =s(x)g.
\end{equation}
Here $\Phi_k \in
\{\BC^p, \, L^2_p(0,l)\}$ ($k=1,2$)
and $\BC$ denotes the complex plane.  This identity
differs from the identity $AS-SA^*=i (\Phi_1 \Phi_2^*+\Phi_2 \Phi_1^*)$
\cite{SaL1, SaL2} for an operator with difference
kernel.   Matrices satisfying a discrete  analogue
of (\ref{0.4}) were treated in \cite{FKS2}.
The operator identity (\ref{0.4}) for
the case, when $k$ in (\ref{0.2})  is a vector, was studied
in \cite{KKL}.
It could be useful also to prove  (\ref{0.4}) directly.
In fact, we prove below a somewhat more general identity. 
\begin {Pn} \label{id}
Let  the operator $S$ in $L^2_p(0,l)$  ($0<l<\infty$)  be defined
by
\begin{equation}\label{v8}
{\cal S}f=f(x)+ \frac{1}{2}\int_0^l
\int^{x+t}_{|x-t|} k\Big(
\frac{r+x-t}{2} \Big)  \wt k \Big( \frac{r+t-x}{2} \Big)d r f(t) d
t,
\end{equation}
where 
$\sup_{0<x<l}\big(\|k(x)\|+\|\wt k(x)\|\big)<\infty$. Then ${\cal S}$ 
satisfies the operator identity
\begin{equation}\label{v9}
A{\cal S}-{\cal S}A^*=i \int_0^l \big(I_p+\psi (x) \wt \psi (t)\big) \, \cdot \, dt,
\end{equation}
where $\psi (x)=\int_0^x k(t)dt,$  $\wt \psi (x)=\int_0^x \wt k(t)dt$.
\end{Pn} 
\begin{proof}.  Using   (\ref{v8})  and changing the order of integration we have
\begin{eqnarray}\label{v10}&&
A{\cal S}f=Af+i \int_0^l \g_1(x,t) f(t)dt, \\
 \label{v11} &&
 \g_1(x,t):= \frac{1}{2} \int_0^x
\int^{y+t}_{|y-t|} k \Big(
\frac{r+y-t}{2} \Big) \wt k \Big( \frac{r+t-y}{2} \Big)d rdy.
\end{eqnarray}
Taking into account that for the scalar product  $(\cdot , \cdot)_l$ in $L^2_p(0,l)$
we have $(A^*f,g)_l$=$(f,Ag)_l$, rewrite ${\cal S}A^*$ in the form
\begin{eqnarray}\label{v12}&&
{\cal S}A^*f=A^*f-i \int_0^l \g_2(x,t) f(t)dt, \\
 \label{bv13} &&
 \g_2(x,t):= \frac{1}{2} \int_0^t
\int^{x+y}_{|x-y|} k\Big(
\frac{r+x-y}{2} \Big)\wt  k \Big( \frac{r+y-x}{2} \Big)d rdy.
\end{eqnarray}
First, consider the case $t \geq x$. From  (\ref{v11}), after changes of variables
$\xi=(r+y-t)/2$ and $\eta=t-y+\xi$, we get
\begin{eqnarray}
 \label{v13} &&
 \g_1(x,t)= \frac{1}{2} \int_0^x
\int^{y+t}_{t-y} k\Big(
\frac{r+y-t}{2} \Big)\wt k \Big( \frac{r+t-y}{2} \Big)d rdy \\
\nonumber &&
=\int_0^x\int_0^y k(\xi)\wt k(t-y+\xi)d \xi dy=\int_0^x\int_{\xi}^x k(\xi) \wt k(t-y+\xi)dy d \xi 
\\ \nonumber &&
=\int_0^x\int^t_{t-x+\xi} k(\xi)\wt k(\eta)d\eta d \xi .
\end{eqnarray}
Next, calculate $\g_2(x,t)$ ($t \geq x$). From  (\ref{bv13})  
it follows that
\begin{equation}\label{v14}
  \g_2(x,t)= \g_{21}(x,t)+ \g_{22}(x,t),
\end{equation}
where
\begin{eqnarray} \label{v15}  &&  \g_{21}(x,t):=
  \frac{1}{2} \int_0^x
\int^{x+y}_{x-y} k\Big(
\frac{r+x-y}{2} \Big)\wt k \Big( \frac{r+y-x}{2} \Big)d rdy,
  \\
 \label{v16} &&
 \g_{22}(x,t):=
  \frac{1}{2} \int_x^t
\int^{x+y}_{y-x} k\Big(
\frac{r+x-y}{2} \Big)\wt k \Big( \frac{r+y-x}{2} \Big)d rdy.
\end{eqnarray}
Replace  the variable $r$ by $\eta =(r+y-x)/2$ in (\ref{v15}) , then change the order
of integration, and after that put  $\xi=x-y+\eta$ and change the order of integration again to obtain
\begin{equation}\label{v17}
\g_{21}(x,t)=
 \int_0^x
\int_{0}^{\xi} 
k(\xi) \wt k(\eta)d\eta d \xi .
\end{equation}
In (\ref{v16}), replace  $r$ by $\xi =(r+x-y)/2$, change the order of integration
and put  $\eta =y-x+\xi$. We get
\begin{equation}\label{v18}
\g_{22}(x,t)=
 \int_0^x
\int_{\xi}^{t-x+\xi} 
k(\xi) \wt k(\eta)d\eta d \xi .
\end{equation}
By   (\ref{v13}),  (\ref{v14}),  (\ref{v17}), and  (\ref{v18})  the equality
 \begin{equation}\label{v19}
\g_{1}(x,t)+ \g_{2}(x,t)=
 \int_0^x
\int_{0}^{t} 
k(\xi)\wt k(\eta)d\eta d \xi =\psi (x) \wt \psi (t)  
\end{equation}  
is true for $t \geq x$. 
Using similar calculations one can show that  (\ref{v19}) holds also for $x \geq t$,
i.e.,  (\ref{v19})  is true for all $0 \leq x,t \leq l$. Finally, formulas  (\ref{v10}),  (\ref{v12}), and  (\ref{v19})
yield (\ref{v9}).
\end{proof}

\section{ISP and semiseparable  operators $S_l$} \label{ISP}
\setcounter{equation}{0}
In this section we consider  matrix functions of the form
\begin{eqnarray}\label{v20}&&
\vp(\la)=i \t_1^*(\la I_n - \b)^{-1} \t_2 e^{- 2 i  \lambda D}R, 
 \\
\label{v20'}&&
D={\mathrm{diag}}
\{d_1,  \ldots, d_p  \}, \quad
 d_{k_1} \geq d_{k_2} \geq 0 \quad {\mathrm{for}} \quad k_1>k_2,
\end{eqnarray}
where $\t_j$ ($j=1,2$) is an $n \times p$ matrix with the $m$-th column denoted by $\t_{j,m}$,
$\b$ is an $n \times n$ matrix,
$R$ is a $p \times p$ matrix, and $D$ is  a $p \times p$ diagonal matrix.
We  do not suppose here that $\t_j$ and $\b$ satisfy the  identity (\ref{v23}). 

\begin{Pn} \label{PEcase}
Let matrix function $\vp$ be given by  (\ref{v20}). Then,  the 
matrix function
$s$, which is defined via $\vp$ by ( \ref{1.3}), has the form $s=CR$,
where $C=\big[c_1 \, \,c_2 \, \ldots \, c_p \big]$,
 the columns $c_m$ $(p \geq m \geq 1)$ being given by the formulas
\begin{eqnarray}\label{v21} &&
c_m(x)= 0  \quad {\mathrm{for}} \quad   0 \leq x  \leq d_m, \\
\label{v21'} &&
c_m(x)=2 \t_1^* \int_0^{x-d_m}\exp \{2it \b \}dt\t_{2,m}  \quad {\mathrm{for}} \quad x \geq d_m , 
\end{eqnarray}  
and the function $\vp$ is the Weyl function of system (\ref{0.1}) with potential $v$ 
satisfying (\ref{v2}).
\end{Pn}
\begin{proof}.
First,
choose $M>0$ such that $\s(\b+ i M I_n)\subset  \BC_+$,
where $\s$ means spectrum and $\BC_+$ is the open upper halfplane. 
According to  (\ref{v20})  $\vp(\la)$ is analytic and the function $\la \vp(\la)$
is bounded in the halfplane $\Im \la< -M$.  So, the conditions   (\ref{v4}) and (\ref{v5})
on $\vp$ are fulfilled.   The fact that $s$ is absolutely
continuous and satisfies conditions (\ref{v6}) and (\ref{v7})  is immediate from
(\ref{v21}) and (\ref{v21'}).
Therefore, after we  have proved   (\ref{v21}) and (\ref{v21'}) , it will  follow from Theorem
\ref{TmssaD}  that  $\vp$ is the Weyl function of system (\ref{0.1}) with potential $v$ 
satisfying (\ref{v2}).

Now, let us prove  (\ref{v21}) and (\ref{v21'}).  As $\la \vp(\la)$
is bounded,  one can rewrite (\ref{1.3}) as a pointwise limit:
 \begin{eqnarray}\label{v22} && s=\big[c_1 \, \,c_2 \, \ldots \, c_p \big]R,  \quad
 c_m(x)= -\frac{1}{\pi} \t_1^*\int_{- \infty}^{\infty}e^{i \la
(x-d_m)}
 \\
 \nonumber &&
\times  \lambda^{-1} (\la I_n -2 \b)^{-1}  d \xi \t_{2,m} \quad (\lambda = \xi +i \eta , \quad \eta <-2M).
\end{eqnarray}
Introduce   the counterclockwise oriented contours, where $\xi$ may take complex values:
\[
\G_{a}^+=[-a,a] \bigcup \{ \xi : \, |\xi |=a, \, \Im \xi >0\}, \, \, \G_{a}^-=[-a,a] \bigcup \{ \xi : \, |\xi |=a,
 \, \Im \xi <0\}.
\]
For $\la=\xi+i \eta$ and for the fixed values of $\eta <-2M$, it follows from  (\ref{v22})  that
\begin{equation}\label{v24} 
c_m(x)= -\frac{1}{ \pi}\theta_{1}^{*}
\lim_{a \to \infty} \int_{\G_a^+ }e^{i \la
(x-d_m)} \lambda^{-1} ( \lambda I_{n} - 2\beta )^{-1}d \xi  \theta_{2,m}
\end{equation}
in the case $x \geq d_m$, and
\begin{equation}
\label{v25} 
c_m(x)= \frac{1}{ \pi}\theta_{1}^{*}
\lim_{a \to \infty} \int_{\G_a^- }e^{i \la
(x-d_m)} \lambda^{-1} ( \lambda I_{n} - 2\beta )^{-1}d \xi  \theta_{2,m}
\end{equation}
in the case $ x \leq d_m$.  As $e^{i \la
(x-d_m)} \lambda^{-1} ( \lambda I_{n} - 2\beta )^{-1}$ is analytic  with
respect to $\xi$  inside $\G_a^-$ and on the contour itself, equality   (\ref{v21}) is immediate from (\ref{v25}).

Next, consider the case  $x \geq d_m$.
For sufficiently large $a$
all the poles of  $( \lambda I_{n} - 2\beta)^{-1}$ (and the pole $\xi=-i \eta$ of $\la^{-1}$)
are contained inside $\G_a^+$ and taking into account (\ref{v24}) we have
\begin{equation}\label{v26} 
c_m(x)= -\frac{1}{ \pi}\theta_{1}^{*}
\int_{\G_a^+ }e^{i \la
(x-d_m)} \lambda^{-1} ( \lambda I_{n} - 2\beta )^{-1}d \xi  \theta_{2,m}.
\end{equation}
Let us approximate $\b$ by matrices $\b_{\ve}$ such that $\|\b -\b_{\ve}\|<\ve$
and $\det \b_{\ve}\not=0$ (if  $\det \b\not=0$ we put $\b = \b_{\ve}$).
It is easy to see that 
\begin{equation}\label{1.14}
\lambda^{-1}( \lambda I_{n} - 2\b_{\ve} )^{-1}=(2\b_{\ve})^{-1}\big(( \lambda I_{n} - 2\b_{\ve} )^{-1}-
  \lambda^{-1}I_n\big).
\end{equation}
For sufficiently small $\ve$
all the poles of  $( \lambda I_{n} - 2\beta_{\ve})^{-1}$ 
are contained inside $\G_a^+$ and we have
\begin{equation}\label{1.15}
\frac{1}{2 \pi i}
\int_{\G_a^+ }e^{i \xi
x} ( \lambda I_{n} - 2\beta_{\ve})^{-1}d \xi=e^{\eta x}\exp({2i x \beta_{\ve}})  \quad
(x \geq 0).
\end{equation}
Finally,  using  (\ref{v26})-(\ref{1.15}) we get 
\[
c_m(x)=\lim_{\ve \to 0}2 \t_1^* \int_0^{x-d_m}\exp \{2it \b_{\ve} \}dt\t_{2,m}.
\]
Hence,
formula (\ref{v21'}) is immediate.
\end{proof}
\begin{Rk}\label{win}
Note that the matrix functions $\vp$ of the form  (\ref{v20})  in general position
do not satisfy  (\ref{v3}) and so they do not satisfy  in a scalar case 
conditions of Lemma 1 \cite{SaA2}, but the conditions of Theorem \ref{TmssaD}
are fulfilled.
\end{Rk}

By Proposition \ref{PEcase} the matrix function $k$ in the expression
(\ref{v27}) for the kernel  of the operator $S_l$, generated by the Weyl
function $\vp$ of the form  (\ref{v20}),  is given by the formula
\begin{eqnarray}\label{v28} &&
k(x)=s^{\prime}(x)=
2\t_1^*e^{2ix \b}\nu \chi(x)R, \quad \nu:=\{ \exp (-2id_m \b )\t_{2,m}\}_{m=1}^p,
\\ \nonumber &&
 \chi(x)={\mathrm{diag}}
\{\chi_1(x), \chi_2(x), \ldots, \chi_k(x) \}, \quad
\chi_m(x):
= \left\{\begin{array}{l}0,  \quad   0 \leq x  < d_m, \\
1, \quad x > d_m .
\end{array} \right.
\end{eqnarray}
According to (\ref{v27}) and (\ref{v28})  we have
\begin{equation}\label{1v2} 
K(x,t)=2\t_1^*\int_{|x-t|}^{x+t}\exp\Big(i(r+x-t)\b\Big)Q(r,x,t)\exp\Big(-i(r+t-x)\b^*\Big)dr\t_1,
\end{equation}
where
\begin{equation}\label{v29} 
Q(r,x,t)=\nu \chi\Big(\frac{r+x-t}{2}\Big)RR^*\chi\Big(\frac{r+t-x}{2}\Big)\nu^*.
\end{equation}
The matrix function  $Q(r,x,t)$
is piecewise constant with respect to $r$ and 
without loss of generality we assume $Q(0,x,t)=0$.
It is easy to see that $Q(r,x,t)$
has only a finite number of  jumps $\{Q_j\}$.
Moreover,
if $\s(\b) \cap \s(\b^*)=\emptyset$, the matrix identity
$i(\b X_j - X_j \b^*)=Q_j$ always has the solution $X_j$.
Therefore we have
\begin{equation}\label{v29'} 
e^{ir \b}Q_je^{-ir \b^*}=\frac{d}{dr}\Big(e^{ir \b}X_je^{-ir \b^*}\Big).
\end{equation}
Hence, according to (\ref{1v2}) and (\ref{v29'}) we can express the kernel
$K(x,t)$ of $S$ explicitly in terms of matrix exponents and $\{X_j\}$. It follows also from
 (\ref{v27})  that
\begin{equation}\label{1v3} 
K(x,t)=K(t,x)^*,
\end{equation}
and so we need to simplify  (\ref{1v2}) only for $x>t$.

Another approach to the presentation of $K$ in terms of matrix exponents is given in the
following lemma.
\begin{La}\label{LaMEx}
Put 
\begin{equation}\label{2v1} 
g_j(r):=[0 \quad e^{i r \b}]e^{r E_j}\left[\begin{array}{c}
I_p\\0\end{array}\right], \quad E_j:=\left[\begin{array}{cc}
-i \b^*&0\\Q_j&-i \b\end{array}\right].
\end{equation}
Then we have
\begin{equation}\label{2v2} 
\Big(\frac{d}{dr}g_j\Big)(r)=e^{ir \b}Q_je^{-ir \b^*}.
\end{equation}
\end{La}
\begin{proof}. By (\ref{2v1}) we have
\[
\frac{d}{dr}g_j=i \b g_j+[0 \quad e^{i r \b}]E_je^{r E_j}\left[\begin{array}{c}
I_p\\0\end{array}\right]=i \b g_j+e^{ir \b}Q_je^{-ir \b^*}-i \b g_j,
\]
and  (\ref{2v2}) is immediate.
\end{proof}
Recall \cite{GGK1} that the operator $S$ is called semiseparable, when $K$ admits
representation 
\begin{equation}\label{2v3} 
K(x,t)=F_1(x)G_1(t) \quad {\mathrm{for}} \, \,x>t, \quad K(x,t)=F_2(x)G_2(t)
 \quad {\mathrm{for}} \, \,x<t,
\end{equation}
where $F_1$ and $F_2$ are $p \times \wt p$ matrix functions and 
$G_1$ and $G_2$ are $\wt p \times  p$ matrix functions.
For the operator $S$ to be semiseparable, assume 
\begin{equation}\label{v31} 
RR^*=I_p.
\end{equation}
Then the matrix function $Q$  has the form
\begin{equation}\label{v32} 
Q(r,x,t)=\nu \chi\Big(\frac{r+t-x}{2}\Big)\nu^* \quad {\mathrm{for}} \quad  x>t.
\end{equation}
Rewrite  (\ref{v20'}) as
\begin{equation}
	\label{1v4}
D={\mathrm{diag}}
\{\wt d_1I_{p_1},  \ldots, \wt d_k I_{p_k}  \}, \quad p_1+ \ldots +p_k=p, \quad
 \wt d_{k_1} > \wt  d_{k_2} \geq 0 \, \, {\mathrm{for}} \, \, k_1>k_2,
\end{equation}
and put
\begin{equation}
	\label{1v5}
Q_j=\nu P_j\nu^*, \quad P_j= {\mathrm{diag}}
\{0,  \ldots, 0, \,  I_{p_j}, \, 0, \ldots, 0  \}.
\end{equation}
\begin{Rk} \label{Rknot} {\rm Some notations.}
Further  we  consider $K(x,t)$  $(x>t)$ on the
intervals $\wt d_m <t<\min(x, \wt d_{m+1})$,
where we choose such $m$ that the inequalities $\wt d_m<x$ hold.
If $\wt d_1>0$ we put $\wt d_0=0$ and include the interval $\wt d_0 <t<\min(x, \wt d_{1})$
into consideration. If  $x>\wt d_k$, we include the interval $\wt d_k<t<\wt d_{k+1}$ $(\wt d_{k+1}=x)$.
Some matrix functions, like $B(t)$ and $C(t)$, will be considered on the intervals as above,
but with $x=l$.
In the following, in all such cases we simply write  $\wt d_m <t<\wt d_{m+1}$.
We also assume that $\sum_{j=1}^m \ldots=0$, when $m=0$.
\end{Rk}
\begin{Pn} \label{D2case}
Let the matrix function $\vp$ be given by  (\ref{v20}), where $D$ satisfies  (\ref{1v4})
and $R$ is unitary.
Assume also that  the matrix identities
\begin{equation}\label{v33} 
i(\b X_j - X_j \b^*)=Q_j ,
\end{equation}
where $Q_j$ are given by  (\ref{1v5}), have solutions $X_j$. Then  the operator $S$,  
which is defined via $\vp$ by formulas (\ref{0.2}), $k=s^{\prime}$ and (\ref{1.3}), is semiseparable,
and its kernel $K(x,t)$ $(0<x,t<l)$ is given by relation
 (\ref{1v3}) and by the equalities
\begin{equation}
 \label{1v7}
K(x,t)=2\t_1^*\left(e^{2ix \b}Z_m e^{-2it \b^*} 
-e^{2i(x-t) \b}\wt Z_m\right)\t_1  \quad  (\wt d_m <t<\wt d_{m+1})
\end{equation}
for $t<x<l$. 
Here
\begin{equation}\label{1v8} 
Z_m=\sum_{j=1}^m X_j, \quad
\wt Z_m=\sum_{j=1}^m\Big(\exp\big(2i\wt d_j\b\big)\Big)X_j\exp\big(-2i\wt d_j\b^*\big).
\end{equation}
Moreover,  there are self-adjoint
solutions of  (\ref{v33}) and we suppose $X_j=X_j^*$ in  (\ref{v33})  and (\ref{1v8}).
\end{Pn}
\begin{proof}. First,  note  that
we can choose $X_j=X_j^*$ because the adjoint of each solution of
(\ref{v33}) also satisfies (\ref{v33}). 

Next, using  (\ref{1v2}),  (\ref{v32}), and   (\ref{1v5}) we get   equalities 
 \begin{eqnarray}\nonumber &&
 K(x,t) = 2\t_1^*\Big(\sum_{j=1}^m\int_{x-t+2\wt d_j}^{x+t}\big(\exp i(r+x-t)\b \big)Q_j 
 \\&&\label{1v9}
\times 
 \exp\big(-i(r+t-x)\b^* \big) dr\Big)\t_1
 \end{eqnarray}
for $t<x$ and $\wt d_m <t<\wt d_{m+1}$.
From  (\ref{v29'}) and  (\ref{1v9})
it follows that   (\ref{1v7}) holds. Formula  (\ref{1v3}) was derived earlier.
\end{proof}
\begin{Rk}
By (\ref{1.6})-(\ref{1.8}) and   (\ref{v28})  the equality
$v(x)=0$ is valid for $ 0<x<d_1$  in the case $d_1>0$. This fact corresponds
to the  inequality
\begin{equation}  \label{1v10}
\sup_{x \leq  d_1, \, \Im \la<-M} \left\|e^{i x \la} e^{i x \la j} 
\left[\begin{array}{c}
\phi(\la)e^{-2i\la D} \\ I_{p}
\end{array}
\right]\right\|<\infty ,
\end{equation}
which can be easily  checked directly and  is implied also by  (\ref{v1}).
\end{Rk}
When the operator
$S=I+\int_0^lK(x,t) \cdot dt$ is semiseparable and its kernel $K$ is given by (\ref{2v3}),
the kernel of the operator $T=S^{-1}$ is expressed in terms
of the $2 \wt p \times 2 \wt p$ solution $U$ of the differential equation
\begin{equation}  \label{2v4}
\Big(\frac{d}{d x}U\Big)(x)=H(x)U(x), \quad x \geq 0, \quad U(0)=I_{2 \wt p}, 
\end{equation}
where
\begin{equation}  \label{2v5}
H(x):=B(x)C(x), \quad B(x)=\left[\begin{array}{c}
-G_1(x) \\ G_2(x)
\end{array}
\right] , \quad
C(x)=\left[\begin{array}{lr}
F_1(x) & F_2(x)
\end{array}
\right].
\end{equation}
Namely, we have (see, for instance, \cite{GGK1})
\begin{equation}  \label{2v6}
T=S^{-1}=I+\int_0^lT(x,t) \, \cdot \, dt, 
\end{equation}
\begin{equation}  \label{2v7}
T(x,t)=\left\{\begin{array}{l}
C(x)U(x)\big(I_{2 \wt p}-P^{\times}\big)U(t)^{-1}B(t), \quad x>t, \\
-C(x)U(x)P^{\times}U(t)^{-1}B(t), \quad x<t. \end{array} \right.
\end{equation}
Here $P^{\times}$ is given in terms of the $\wt p \times \wt p$ blocks
$U_{21}(l)$ and $U_{22}(l)$ of $U(l)$:
\begin{equation}  \label{2v8}
P^{\times}=\left[\begin{array}{lr}
0 & 0 \\ U_{22}(l)^{-1}U_{21}(l) & I_{\wt p}
\end{array}
\right],
\end{equation}
and the invertibility of  $U_{22}(l)$ is a necessary and sufficient condition
for the invertibility of $S$.  

If $K$ admits  the representation
\begin{equation}  \label{2v9}
K(x,t)=\left\{\begin{array}{l}
Ce^{x{\cal A}} \big(I_{2\wt p}-P\big)e^{-t{\cal A}} B, \quad x>t, \\
-Ce^{x{\cal A}} Pe^{-t{\cal A}} B, \quad x<t, \end{array} \right.
\end{equation}
where ${\cal A}$, $B$, and  $C$ are constant matrices, then $U$
is calculated explicitly \cite{GK84}. In our case a representation
\begin{equation}  \label{2v10}
K(x,t)=\left\{\begin{array}{l}
C_me^{x{\cal A}} \big(I_{2\wt p}-P_m\big)e^{-t{\cal A}} B_m, 
\quad t<x<l, \,  \wt d_m <t<\wt d_{m+1},\\
-C_me^{x{\cal A}} P_me^{-t{\cal A}} B_m, 
\quad x<t<l, \,  \wt d_m <x<\wt d_{m+1},
\end{array} \right.
\end{equation}
where $\wt p =n$ and
\begin{equation}  \label{2v11}
{\cal A}=2i \left[\begin{array}{lr}
\b & 0 \\ 0 & \b^*
\end{array}
\right],
\end{equation}
easily follows from  (\ref{1v3})  and   (\ref{1v7}).
However, (\ref{2v10}) is insufficient for the explicit construction of $U$
and we shall construct $U$ and $T$ explicitly, using more general formulas
(\ref{2v4})-(\ref{2v8}). For this purpose we introduce $B(x)$ and $C(x)$ 
($0<x<l$)
by the equalities
\begin{equation}  \label{2v12}
B(x)=\sqrt{2}\left[\begin{array}{c}e^{-2ix \b}\wt Z_m-Z_me^{-2ix \b^*} \\
e^{-2ix \b^*}
\end{array}
\right]\t_1  \quad (\wt d_m <x<\wt d_{m+1}),
\end{equation}
\begin{equation}  \label{2v13}
C(x)=\sqrt{2}\t_1^*
\left[\begin{array}{lr}e^{2ix \b} &
e^{2ix \b} Z_m- \wt Z_m e^{2ix \b^*} \end{array}
\right]  \quad (\wt d_m <x<\wt d_{m+1}),
\end{equation}
where $Z_m=Z_m^*$ and $\wt Z_m=\wt Z_m^*$ are defined in (\ref{1v8}).

\begin{Pn} \label{PnInv} Let the conditions of Proposition \ref{D2case}
be fulfilled and let $S$ be defined via $\vp$ by formulas (\ref{0.2}), $k=s^{\prime}$ and (\ref{1.3}).
Then the operator $T=S^{-1}$ is given by formulas
(\ref{2v6})-(\ref{2v8}), (\ref{2v12}), (\ref{2v13}), and
\begin{equation}  \label{2v14}
U(x)=\Om_m e^{-x{\cal A}}e^{x{\cal A}_m^{\times}}\Xi_m^{-1}U(\wt d_m)
\quad  (\wt d_m \leq x \leq \wt d_{m+1}), \quad U(0)=I_{2n},
\end{equation}
where ${\cal A}$ is defined by (\ref{2v11}) and
\begin{equation}  \label{2v15}
{\cal A}_m^{\times}:={\cal A}+2Y_m, \quad Y_m:=
\left[\begin{array}{c}\wt Z_m\\
 I_n
\end{array}
\right] \t_1\t_1^*
\left[\begin{array}{lr}I_n & -\wt Z_m
\end{array}
\right],
\end{equation}
\begin{equation}  \label{2v16}
\Om_m:=\left[\begin{array}{lr}I_n & -Z_m\\
0 & I_n
\end{array}
\right],   \quad 
 \Xi_m:=\Om_m e^{-\wt d_m{\cal A}}e^{\wt d_m {\cal A}_m^{\times}}.
\end{equation}
Moreover, we have
\begin{equation}  \label{2v17}
U(x)^*JU(x)=J, \quad U(x)^{-1}=JU(x)^*J^*, \quad
J:=\left[\begin{array}{lr}0 &-I_n \\
I_n & 0
\end{array}
\right].
\end{equation}
\end{Pn}
\begin{proof}. Recall that $B$ and $C$ are recovered from $K$
by formulas  (\ref{2v3}) and (\ref{2v5}). In view of (\ref{1v3}) and (\ref{1v7}) 
we have
\begin{equation}  \label{2v18}
F_1(x)= \sqrt{2}\t_1^*
e^{2ix \b} , \qquad F_2(x)= G_1(x)^*, 
\end{equation}
\begin{eqnarray}  \label{2v19}
G_1(x)&=& \sqrt{2}\big( Z_me^{-2ix \b^*}-e^{-2ix \b}\wt Z_m\big)\t_1
\quad  (\wt d_m < x < \wt d_{m+1}), \\  \label{2v20}
\quad G_2(x)&=&F_1(x)^*.
\end{eqnarray} 
Therefore, formulas   (\ref{2v5}) and  (\ref{2v18})-(\ref{2v20})  imply
that $B$ and $C$ corresponding to $S$ are given by
 (\ref{2v12}) and (\ref{2v13}). It follows from  (\ref{2v5}) and
 (\ref{2v11})-(\ref{2v13}) that
 \begin{equation}  \label{2v21}
H(x)=2\Om_m e^{-x{\cal A}}Y_me^{x{\cal A}}\Om_m^{-1} \quad  (\wt d_m < x < \wt d_{m+1}), 
\end{equation}
where $Y_m$ is given in (\ref{2v15}), $\Om_m$ is given in (\ref{2v16}),
and
\begin{equation}  \label{2v22}
\Om_m^{-1}=\left[\begin{array}{lr}I_n & Z_m\\
0 & I_n
\end{array}
\right].
\end{equation}
According to (\ref{2v14}), (\ref{2v15}), and (\ref{2v21}) we get
\[
\Big(\frac{d}{d x}U\Big)(x)=\Om_m e^{-x{\cal A}}\Big({\cal A}_m^{\times}-{\cal A}\Big)
e^{x{\cal A}_m^{\times}}\Xi_m^{-1}U(\wt d_m)=H(x)U(x) 
\]
for $\wt d_m < x < \wt d_{m+1}$, and so $U$ of the form  (\ref{2v14}) satisfies (\ref{2v4}).
In other words,  formulas (\ref{2v12})-(\ref{2v14}) define explicitly $B$, $C$ and $U$,
which are used in the expressions   (\ref{2v7}) and (\ref{2v8}) to construct the kernel
of $T=S^{-1}$.

It remains to prove  (\ref{2v17}). Note that
\begin{equation}  \label{2v23}
J {\cal A}^*J^*=-{\cal A}, \quad J\Om_m^*J^*=\Om_m^{-1}, \quad JY_m^*J^*=-Y_m.
\end{equation}
Hence, we have $JH^*J^*=-H$, i.e., 
\begin{equation}  \label{2v24}
\frac{d}{dx}\big(U(x)^*JU(x)\big)\equiv 0.
\end{equation}
Formula   (\ref{2v17}) follows from (\ref{2v24}) and from $U(0)=I_{2n}$.
\end{proof}

Taking into account  (\ref{2v12})-(\ref{2v14}) we get
\begin{eqnarray} && \label{2v25}
\wt F(x):=C(x)U(x)=\sqrt{2}\t_1^*[I_n \quad - \wt Z_m]e^{x{\cal A}_m^{\times}}\Xi_m^{-1}U(\wt d_m),
\\ && \label{2v26}
\wt G(t):=U(t)^{-1}B(t)=\sqrt{2}U(\wt d_m)^{-1}\Xi_m e^{-t{\cal A}_m^{\times}}
\left[\begin{array}{c}\wt Z_m\\
 I_n
\end{array}
\right] \t_1,
\\
&& \nonumber \wt d_m < x < \wt d_{m+1}, \quad \wt d_m < t< \wt d_{m+1}.
\end{eqnarray}
\begin{Cy} Let the conditions of Proposition \ref{D2case}
be fulfilled. Then the kernel $T(x,t)$ of the operator $T=S_l^{-1}$  has the form
\begin{equation}  \label{2v27}
T(x,t)=\left\{\begin{array}{l}\wt F(x)
\big(I_{2  n}-P^{\times}\big)\wt G(t), \quad x>t, \\
-\wt F(x)P^{\times}\wt G(t), \quad x<t,\end{array} \right.
\end{equation}
where $\wt F$ and $\wt G$ are given by (\ref{2v25}) and (\ref{2v26}).
\end{Cy}
By (\ref{v28}), (\ref{2v15}), and (\ref{2v26}) for $\wt d_m < t< \wt d_{m+1}$ we get
\begin{eqnarray}\nonumber
\wt G(t) k(t)&=&\sqrt{2}U(\wt d_m)^{-1}\Xi_m e^{-t{\cal A}_m^{\times}}(2Y_m)e^{t{\cal A}}\left[\begin{array}{c} I_n\\
 0
\end{array}
\right]\nu \sum_{j=1}^m P_j R \\ \label{3v8} &=&- \sqrt{2}U(\wt d_m)^{-1}
\Xi_m \frac{d}{dt}\Big(e^{-t{\cal A}_m^{\times}}e^{t{\cal A}}\Big)\left[\begin{array}{c} I_n\\
 0
\end{array}
\right]\nu \sum_{j=1}^m P_j R .
\end{eqnarray}
 From Remark \ref{RkV} and formulas (\ref{2v6}), (\ref{2v27}),
and (\ref{3v8})  the explicit solution of the inverse problem is immediate.
\begin{Tm} \label{main}
Let the Weyl matrix function $\vp$ be given by
 (\ref{v20}), where $D$ satisfies  (\ref{1v4})
and $R$ is unitary.
Assume also that  the matrix identities (\ref{v33}),
where $Q_j$ are given by  (\ref{1v5}), have solutions $X_j=X_j^*$.
Then the ISP solution $v$ is given by the formula
\begin{eqnarray}&&\label{3v9} 
v(l)=k(l)+\wt F(l)\big(I_{2n}-P^{\times}\big)\sum_{m=1}^{N} \sqrt{2}U(\wh d_m)^{-1}\Xi_m
\\  && \nonumber
\times\Big(
e^{-\wh d_m{\cal A}_m^{\times}}e^{\wh d_m{\cal A}}-
e^{-\wh d_{m+1}{\cal A}_m^{\times}}e^{\wh d_{m+1}{\cal A}}\Big)\left[\begin{array}{c} I_n\\
 0
\end{array}
\right]\nu \sum_{j=1}^m P_j R ,
\end{eqnarray}
where $k$ is given by (\ref{v28}), $U$ is given by (\ref{2v14}),  $P^{\times}$ is given by
(\ref{2v8}), and $\Xi_m$  is given by (\ref{2v16}). The number $N$ in the sum
is chosen in the following way: if $l<\wt d_1$  then $N=0$; if  $\wt d_j<l<\wt d_{j+1}$ then $N=j$; if $l>\wt d_k$ then
$N=k$. We put $\wh d_m=\wt d_m$ for $m\leq  N$ and $\wh d_{N+1}=l$.
\end{Tm}

{\bf Acknowledgement.}
The work of A.L. Sakhnovich was supported by the Austrian Science Fund (FWF) under
Grant  no. Y330.

\begin{flushright} \it
B. Fritzsche,  \\ 
Fakult\"at f\"ur Mathematik und Informatik, \\
Mathematisches Institut, Universit\"at Leipzig, \\ 
Johannisgasse 26,  D-04103 Leipzig, Germany,\\ 
e-mail: {\tt fritzsche@math.uni-leipzig.de } \\  $ $ \\

B. Kirstein, \\
Fakult\"at f\"ur Mathematik und
Informatik, \\
Mathematisches Institut, Universit\"at Leipzig,
\\ Johannisgasse 26,  D-04103 Leipzig, Germany, \\ 
e-mail: {\tt kirstein@math.uni-leipzig.de } \\  $ $ \\

A.L. Sakhnovich, \\  Fakult\"at f\"ur Mathematik,
Universit\"at Wien,
\\
Nordbergstrasse 15, A-1090 Wien, Austria \\ 
e-mail: {\tt al$_-$sakhnov@yahoo.com }
\end{flushright}

\end{document}